# A SIMPLE METHOD FOR SEARCHING FOR PRIME PAIRS IN THE GOLDBACH CONJECTURE

## WEI SHENG ZENG AND ZIQI SUN

ABSTRACT. In this paper we introduce a simple method of searching for the prime pairs in the famous Goldbach Conjecture. The method, which is based on certain integer identities as well as an observation related to the remainder property, enables us to cut the search length by half and therefore can be considered as a short cut in searching for the prime pairs in the Goldbach Conjecture.

## 1. INTRODUCTION

The Goldbach Conjecture is one of the oldest and best-known unsolved problems in number theory and in all of mathematics. It states: Every even integer greater than 2 can be expressed as the sum of two primes. In this paper we introduce a simple method of finding prime pairs conjectured in this conjecture.

Let $N$ be an even integer. Traditionally, when one wants to find the prime pair $\{p_1, p_2\}$ with $N = p_1 + p_2$ in the Goldbach Conjecture, one needs to search for such a pair from all the integer pairs of the form $\{p, N - p\}$, where $p$ is the prime less than $N$. Based on the prime-counting function, the number of primes less than or equal to $N$ is approximately given by $N/lnN$, as $N$ goes to the infinity. Thus, roughly speaking, this lengthy process involves roughly $N/2lnN$ possible integer pairs to be tested. In this paper, we propose a new procedure of finding the prime pair in the Goldbach Conjecture in which the number of integer pairs to be tested will be no more than half of those pairs to be tested in the traditional approach.

The argument developed in this paper is based on certain integer identities as well as an observation related to the remainder property, which we will establish in Section 2. In Section 3 we shall relate the argument developed in Section 2 to the problem of searching for prime pairs in the Goldbach Conjecture. After we reformulate the argument in a slightly simplified form in Section 4, we will propose a new procedure of searching for the prime pairs in the Goldbach Conjecture in Section 5. Finally, in Section 6, we will test our method on a series of even integers.

## 2. BASIC INTEGER IDENTITIES AND REPRESENTATION FORMULAS

We first arrange all integers ( $\geq 2$) in the following table as shown. The first column lists the row numbers and all the natural numbers are arranged in Column A to Column F, where $y = x - 1$.

| $x$ | A | B | C | D | E | F |
|---|---|---|---|---|---|---|
| 1 | 2 | 3 | 4 | 5 | 6 | 7 |
| 2 | 8 | 9 | 10 | 11 | 12 | 13 |
| 3 | 14 | 15 | 16 | 17 | 18 | 19 |
| 4 | 20 | 21 | 22 | 23 | 24 | 25 |
| 5 | 26 | 27 | 28 | 29 | 30 | 31 |
| 6 | 32 | 33 | 34 | 35 | 36 | 37 |
| 7 | 38 | 39 | 40 | 41 | 42 | 43 |
| .. | .. | .. | .. | .. | .. | .. |
| $x$ | $2+6y$ | $3+6y$ | $4+6y$ | $5+6y$ | $6+6y$ | $7+6y$ |
| .. | .. | .. | .. | .. | .. | .. |

The arrangement of the integers in the above table obeys the following rules:

R1     In each column, integers increase by 6 from row to row.
R2     In each row, integers increase by 1 from column to column.

We observe that

P1     The integers in Column A, C, and E are all even.
P2     In Column B every integer is a multiple of 3. So the integers in Column B are all composite numbers except the first integer 3.
P3     All the primes are contained in Column D and F.

We also observe the following relation between the columns and the remainders when the integers are divided by 6:

| Column | A | B | C | D | E | F |
|---|---|---|---|---|---|---|
| Remainder | 2 | 3 | 4 | 5 | 0 | 1 |

Based on the above observation we shall first explore some relations behind the columns A, B, C, E, D, and F. In the following we shall use the capital letters $A_m$, $B_m$, $C_m$, $D_m$, $E_m$, and $F_m$ to express integer m in the corresponding columns A, B, C, D, E, and F. For example, $D_5$ represents the integer 5 since 5 is in Column D, and $C_{52}$ represents the integer 40 since 40 is in Column C. Similarly, the equation $A_{2+6(x-1)} = 3 + D_{5+6(x-2)}$ actually means the integer $2 + 6(x - 1)$, which is in Column A, equals to the integer 3 plus the integer $(5 + 6(x - 2))$, which is in Column D. So the real equation is $2 + 6(x - 1) = 3 + (5 + 6(x - 2))$, which is just an identity.

Given an even integer $N$ ($\geq 14$) in Column A, then $N$ has a remainder 2 when $N$ is divided by 6. So we can write $N = A_{2+6(x-1)}$ for some integer $x \geq 3$. It is easily verified that the following identity holds for any $n$ and $x$:

$$2 + 6(x - 1) = 3 + 5 + 6(x - 2) = 7 + 6(n - 1) + 7 + 6(x - 2 - n),$$

where $n$ is an integer to be specified later. So we conclude that the even integer $A_{2+6(x-1)}$ can be expressed either by $3 + D_{5+6(x-2)}$ or by $F_{7+6(n-1)} + F_{7+6(x-2-n)}$. Thus we obtain the following representation formula:

(2.1) $\qquad A_{2+6(x-1)} = 3 + D_{5+6(x-2)}$ or $F_{7+6(n-1)} + F_{7+6(x-2-n)}$

One can also determine the range of $n$. Letting $7 + 6(n-1) = 7 + 6(x-2-n)$ we obtain $n = (x-1)/2$. So when $x$ is even, $1 \leq n \leq (x-2)/2$ and when $x$ is odd, $1 \leq n \leq (x-1)/2$. This means that if we want to represent the even integer $A_{2+6(x-1)}$ as a sum of two integers in Column F, there are $n$ possible ways to do that, and $n$ is determined by the above inequalities.

Similarly in Column C, an even integer $N$ has a remainder 4 when N is divided by 6. So we can write $N = C_{4+6(x-1)}$ for $N \geq 16$, where $x \geq 3$. It is easily verified the following identity holds for any n and x:

$$4 + 6(x-1) = 3 + 7 + 6(x-2) = 5 + 6(n-1) + 5 + 6(x-1-n)$$

So we conclude that the even integer $C_{4+6(x-1)}$ can be expressed either by $3 + F_{7+6(x-2)}$ or by $D_{5+6(n-1)} + D_{5+6(x-1-n)}$. Thus we obtain the following representation formula:

(2.2) $\qquad C_{4+6(x-1)} = 3 + F_{7+6(x-2)}$ or $D_{5+6(n-1)} + D_{5+6(x-1-n)}$.

By letting $5 + 6(n-1) = 5 + 6(x-1-n)$ we obtain $n = x/2$, so we can get the range of $n$: When $x$ is odd, $1 \leq n \leq (x-1)/2$ and when $x$ is even, $1 \leq n \leq x/2$. So we obtain the number of ways to express the even integer $C_{4+6(x-1)}$ as a sum of two integers in Column D.

For Column E, an even integer $N$ has a remainder 0 when $N$ is divided by 6. So we can write $N = E_{6+6(x-1)}$ for $N \geq 18$, where $x \geq 3$. It is easily verified the following identity holds for any n and $x$:

$$6 + 6(x-1) = 5 + 6(n-1) + 7 + 6(x-1-n)$$

$$= 7 + 6(n-1) + 5 + 6(x-1-n)$$

So we conclude that the even integer $E_{6+6(x-1)}$ can be expressed either by $D_{5+6(n-1)} + F_{7+6(x-1-n)}$ or by $F_{7+6(n-1)} + D_{5+6(x-1-n)}$. Thus we obtain the following formula:

(2.3) $\qquad E_{6+6(x-1)} = D_{5+6(n-1)} + F_{7+6(x-1-n)}$ or $F_{7+6(n-1)} + D_{5+6(x-1-n)}$.

In this case the range of $n$ is quite clear: $1 \leq n \leq x-1$. This gives number of ways to express the even integer $E_{6+6(x-1)}$ as a sum of an integer in Column D and another integer in Column F.

# 3. REPRESENTATIONS OF EVEN INTEGERS BY SUM OF PRIME PAIRS

Using the representation formulas in $(2.1)$-$(2.3)$, we can now to represent an even integer by the sum of two primes. We can show that if such a representation exists as conjectured by the Goldbach Conjecture, it will be given by one of the representation formulas $(2.1)$-$(2.3)$.

**Theorem 1**. Let N be an even integer greater than or equal to $14$ and if there are primes $p_1$ and $p_2$ so that $N = p_1 + p_2$, then such prime representation is determined by the formulas $(2.1)$-$(2.3)$. Namely,

Case 1. If $N$ has a remainder 2 when it is divided by 6, then either $p_1$ or $p_2$ is 3 and the other is in Column D, or both $p_1$ and $p_2$ are in Column F.

Case 2. If $N$ has a remainder 4 when it is divided by 6, then either $p_1$ or $p_2$ is 3 and the other is in Column F, or both $p_1$ and $p_2$ are in Column D.

Case 3. If N has a remainder 0 when it is divided by 6, then either $p_1$ or $p_2$ is in Column D and the other is in Column F.

Proof. The possibility of such a representation is already shown in the formulas $(2.1)$-$(2.3)$. We only need to prove that beyond such representations there is no other way to represent an even integer by the sum of two primes. We argue by three cases.

In the case 1, the even integer $N = A_{2+6(x-1)}$ has a reminder 2. If it is expressed by the sum of two primes $p_1 + p_2$, then the remainders of the primes $p_1$ and $p_2$ must be either 3 and 5 or 1 and 1. This indicates that $N$ can only be expressed by a prime in Column B, which is 3, and a prime in Column D, or by two primes in Column F. Therefore in this case, the representation is given by the formula $(2,1)$.

In the case 2, the even integer $N = C_{4+6(x-1)}$ has a reminder 4. If it is expressed by the sum of two primes $p_1 + p_2$, then the remainders of the primes $p_1$ and $p_2$ must be either 3 and 1 or 5 and 5. This indicates that N can only be expressed by a prime in Column B, which is 3, and a prime in Column F, or by two primes in Column D. Therefore in this case, the representation is given by the formula $(2.2)$.

In the case 3, the even number $N = E_{6+6(x-1)}$ has a reminder 0. If it is expressed by the sum of two primes $p_1 + p_2$, then the remainder of the primes $p_1$ and $p_2$ must be either 2 and 4, or 3 and 3. Note that the latter case only works for $6 = 3 + 3$, so there is only one possibility: the remainder of one of the primes $p_1$ and $p_2$ is 2 and the remainder of the other one is 4. This indicates that $N$ can only be expressed by a prime in Column D and a prime in Column F. Therefore in this case, the representation is given by the formula $(2.3)$. □

# 4. THE SIMPLIFIED REPRESENTATION FORMULAS

In order to use the representation formulas established in Section 1 and 2 to search the prime pairs conjectured by the Goldbach Conjecture, we need to simplify the formula $(2.1)$-$(2.3)$ by eliminating $x$ in the expressions. In each formula, we can express $x$ in term of $N$ and make a substitution on the right side of the formula. We now carry out this procedure in details.

Case 1. If N is an even integer in Column A (with a remainder 2), then from the formula $(2.1)$ we have $N = 2 + 6(x-1)$ and therefore $x = (N+4)/6$. Thus,

$$F_{7+6(n-1)} + F_{7+6(x-2-n)} = F_{6n+1} + F_{7+6((N+4)/6-2-n)} = F_{6n+1} + F_{N-6n-1}$$

and also

$$3 + D_{5+6(x-2)} = 3 + D_{N-3}.$$

We now have that

(4.1) $$N = 3 + D_{N-3} \text{ or } F_{6n+1} + F_{N-6n-1},$$

where the range of $n$ is determined in the formula $(2.1)$.

Case 2. If $N$ is an even integer in Column C (with a remainder 4), then from the formula $(2.2)$ we have $N = 4 + 6(x-1)$ and therefore $x = (N+2)/6$. Thus,

$$D_{5+6(n-1)} + D_{5+6(x-1-n)} = D_{6n-1} + D_{5+6((N+2)/6-1-n)} = F_{6n-1} + F_{N-6n+1}$$

and also

$$3 + F_{7+6(x-2)} = 3 + F_{N-3}.$$

We now have that

(4.2) $$N = 3 + F_{N-3} \text{ or } D_{6n-1} + D_{N-6n+1},$$

where the range of n is determined in the formula $(2.2)$.

Case 3. If N is an even integer in Column C (with a remainder 0), then from the formula $(2.3)$ we have $N = 6 + 6(x-1)$ and therefore $x = N/6$. Thus,

$$D_{5+6(n-1)} + F_{7+6(x-1-n)} = D_{6n-1} + F_{7+6(N/6-1-n)} = D_{6n-1} + F_{N-6n+1}$$

and

$$F_{7+6(n-1)} + D_{5+6(x-1-n)} = F_{6n+1} + D_{5+6(N/6-1-n)} = F_{6n+1} + D_{N-6n-1}$$

We now have that

(4.3) $$N = D_{6n-1} + F_{N-6n+1} \text{ or } F_{6n+1} + D_{N-6n-1},$$

where the range of $n$ is determined in the formula $(2.3)$.

We summarize the above result in the following theorem.

**Theorem 2**. Let $N$ be an even integer greater than or equal to $14$ and if the Goldbach Conjecture holds for $N$, then $N$ admits one of the representations given by $(4.1) - (4.3)$.

## 5. SEARCHING FOR PRIME PAIRS IN THE GOLDBACH CONJECTURE

Based on the above analysis, we now give a procedure of finding prime pairs in the Goldbach Conjecture:

**Step 1** For a given even integer $N$, find its remainder when it is divided by 6. The remainder should be $2, 4$ or $0$.
**Step 2** Determine the column where $N$ belongs to. The column should be A, C, or E.
**Step 3** Use the formula $(4.1) - (4.3)$ to find the prime pairs.

We now discuss in detail the procedure again in three cases.

Case 1. The even integer $N$ is in Column A, i.e., then $N$ has a remainder $2$ when it is divided by 6. By the formula $(4.1)$, $N = 3 + D_{N-3}$. If $N - 3$ is a prime, then $\{3, N - 3\}$ is the desired prime pair. If $N - 3$ is not a prime, then by the formula $(4.1)$ again,

$$N = F_{6n+1} + F_{N-6n-1}.$$

This means the desired prime pair must be in the form $\{6n + 1, N - 6n - 1\}$ with $n$ chosen in the range given in $(2.1)$. One can test these integer pairs to get a prime pair.

Case 2. The even integer $N$ is in Column C, i.e., then $N$ has a remainder $4$ when it is divided by 6. By the formula $(4.2)$, $N = 3 + F_{N-3}$. If $N - 3$ is a prime, then $\{3, N - 3\}$ is the desired prime pair. If $N - 3$ is not a prime, then by the formula $(4.2)$ again,

$$N = D_{6n-1} + D_{N-6n+1}.$$

This means the desired the prime pair must be in the form $\{6n - 1, N - 6n + 1\}$ with $n$ chosen in the range given in $(2.2)$. One can test these integer pairs to get a prime pair.

Case 3. The even integer $N$ is in Column E, i.e., then $N$ has a remainder $0$ when it is divided by 6. By the formula $(4.3)$,

$$N = D_{6n-1} + F_{N-6n+1} \text{ or } F_{6n+1} + D_{N-6n-1}.$$

This means the desired the prime pair must be in the form $\{6n - 1, N - 6n + 1\}$ or $\{6n + 1, N - 6n - 1\}$ with $n$ chosen in the range given in $(2.3)$. One can test these integer pairs to get a prime pair.

The above procedure can be considered as a short cut in finding the prime pairs in the Goldbach Conjecture. As we discussed earlier in the beginning of the paper, the traditional way of finding such a pairs requires checking on the integer pairs of the form $\{p, N - p\}$, where $p$ is the prime less than N, and based on the prime-counting function, so there are up to approximately $N/2lnN$ such pairs to be checked. In the new procedure given above, since one only focus on one column of primes, so the number of integer pairs to be tested will be reduced by half.

The argument used in this paper is based on the remainder property when an even number is divided by 6. We call this number, 6, is the base number. It remains to see whether one can use other base number to produce better result in shortening the length of searching for the prime pairs in the Goldbach Gonjecture. Our experience indicates that a good base number should be in the form of $2p_1p_2...p_m$ where $p_1$, $p_2$, ... $p_m$ are distinct primes. For instance, the next base number to try is $30 = 2 \cdot 3 \cdot 5$. We shall discuss the further result in this direction in another paper.

## 6. COMPUTATIONAL RESULTS

We now demonstrate the procedure proposed in Section 4 by finding the prime pairs for some even numbers.

From the design of the procedure, we know that the method is for even integers greater than 12. However, the existence of prime pairs is obvious for the even integers $4, 6, 8, 10$, and $12$.

We now come to work on even numbers that are greater than 12. In the process, we use bold numbers to denote composite numbers and unframed number to denote primes.

$14 = 3 + 11$
$16 = 3 + 13$
$18 = 13 + 5 = 11 + 7$
$20 = 3 + 17 = 7 + 13$
$22 = 3 + 19 = 5 + 17$
$24 = 5 + 19 + 7 + 17$
$26 = 3 + 23 = 7 + 19$
$28 = 3 + \mathbf{25} = 11 + 17 = 5 + 23$
$30 = 5 + \mathbf{25} = 11 + 19 = 7 + 23 = 13 + 17$
$32 = 3 + 29 = 7 + \mathbf{25} = 13 + 19$
$34 = 3 + 31 = 5 + 29 = 11 + 23$
$36 = 3 + 31 = 11 + \mathbf{25} = 17 + 19 = 7 + 29 = 13 + 23$
……

When $N = 278$, the remainder of $N$ is 2, so 278 is in Column A :

$N = 3 + D_{N-3} = 3 + 275$

$N = F_{6n+1} + F_{N-6n-1} = 7 + 271 = 13 + \mathbf{265} = 19 + \mathbf{259} = 25 + \mathbf{253} = 31 + \mathbf{247}$
$= 37 + 241 = 43 + \mathbf{235} = \mathbf{49} + 229 = \mathbf{55} + 223 = 61 + \mathbf{217} = 67 + 211 = 73 + \mathbf{205}$
$= 79 + 199 = \mathbf{85} + 193 = \mathbf{91} + \mathbf{187} = 97 + 181 = 103 + \mathbf{175} = 109 + \mathbf{169} = \mathbf{115} + 163 = \mathbf{121} + 157 = 127 + 151 = \mathbf{133} + \mathbf{145}$.

Similarly, when $N = 280$, the remainder of $N$ is 4, so 280 is in Column C:

$N = 3 + C_{N-3} = 3 + 277$
$N = D_{6n-1} + D_{N-6n+1} = 5 + \mathbf{275} = 11 + 269 = 17 + 263 = 23 + 257 = 29 + 251$
$= \mathbf{35} + 245 = 41 + 239 = 47 + 233 = 53 + 227 = 59 + \mathbf{221} = 65 + 215 = 71 + \mathbf{209} = \mathbf{77} + \mathbf{203} = 83 + 197 = 89 + 191 = \mathbf{95} + \mathbf{185} = 101 + 179 = 107 + 173 = 113 + 167 = \mathbf{119} + \mathbf{161} = \mathbf{125} + \mathbf{155} = 131 + 149 = 137 + \mathbf{143}$.

Again, when $N=282$, the remainder of $N$ is 0, so 282 is in Column E:

$N = D_{6n-1} + F_{N-6n+1} = 5 + 277 = 11 + 271 = 17 + \mathbf{265} = 23 + \mathbf{259} = 29 + \mathbf{253}$
$= \mathbf{35} + \mathbf{247} = 41 + 241 = 47 + \mathbf{235} = 53 + 229 = 59 + 223 = \mathbf{65} + \mathbf{217} = 71 + 211 = \mathbf{77} + \mathbf{205} = 83 + 199 = 89 + 193 = \mathbf{95} + \mathbf{187} = 101 + 181 = 107 + \mathbf{175} = 113 + \mathbf{169} = 119 + 163 = \mathbf{125} + 157 = 131 + 151 = 137 + \mathbf{145}$.
$N = F_{6n+1} + D_{N-6n-1} = 7 + \mathbf{275} = 13 + 269 = 19 + 263 = \mathbf{25} + 257 = 31 + 251$
$= 37 + \mathbf{245} = 43 + 239 = \mathbf{49} + 233 = \mathbf{55} + \mathbf{227} = 61 + \mathbf{221} = 67 + 215 = 73 + \mathbf{209} = 79 + \mathbf{203} = \mathbf{85} + 197 = \mathbf{91} + 191 = 97 + \mathbf{185} = 103 + 179 = 109 + 173 = \mathbf{115} + 167 = \mathbf{121} + \mathbf{161} = 127 + \mathbf{155} = \mathbf{133} + 149 = 139 + \mathbf{143}$.

CHONGQING BROADCAST TECHNOLOGY CENTER, CHONGQING, CHINA
  *E-mail address*: zenweyshing@sina.com

DEPARTMENT OF MATHEMATICS, STATISTICS, AND PHYSICS, WICHITA STATE UNIVERSITY, WICHITA, KANSAS, USA
  *E-mail address*: ziqi.sun@wichita.edu